\documentclass[a4paper,12pt]{article}
\input epsf
\newcommand\implies{ \Rightarrow }

\title{Bifurcation in two-dimensional fixed point subspaces}

\author{P.C. Matthews,\\
School of Mathematical Sciences, \\
University of Nottingham, Nottingham NG7 2RD, UK \\
\small paul.matthews@nottingham.ac.uk
}

\begin{document}

\maketitle

\begin{abstract}
Bifurcation with symmetry is considered in the case of an isotropy
subgroup with a two-dimensional fixed point subspace 
and non-zero quadratic terms. 
In general, there are one or three branches of solutions,
and five qualitatively different phase portraits,
provided that two non-degeneracy conditions are satisfied.
Conditions are also derived to determine which of the 
five possible phase portraits occurs, given the coefficients of the quadratic
terms. The results are applied to the problem of bifurcation with
spherical symmetry, where there are six irreducible representations for which
the subspace of solutions with cubic symmetry is two-dimensional.
In each case, the number of solutions and their stability is found.

\end{abstract}

\section{Introduction}

The topic of bifurcation with symmetry 
has wide-ranging applications, including buckling of rods, convection
patterns in fluids and  the structure of tumours and embryos.
In the physical context, a highly symmetric state 
becomes unstable as a control parameter is varied, leading to a state
with reduced symmetry. If there is a high degree of symmetry, 
then the stability problem is degenerate, in the sense that 
several eigenvalues pass through zero simultaneously as the parameter
is varied; the number of eigenvalues is equal to the dimension of one
of the irreducible representations of the original symmetry group. 
This means that the dynamics in the neighbourhood of the bifurcation
is governed by a number of coupled, nonlinear equations 
which in general cannot be solved. 

This complicated multi-dimensional problem can be simplified by
restricting attention to 
the flow-invariant fixed-point subspace of a subgroup of the original 
symmetry group. If this subspace is one-dimensional, 
then the problem reduces to a single differential equation,
so that the usual methods of bifurcation theory can be applied,
leading to a bifurcation of transcritical or pitchfork type. 
This is essentially the equivariant branching lemma
(Vanderbauwhede 1980, Cicogna 1981, Golubitsky,  Stewart and
Schaeffer 1988): if the fixed-point subspace is one-dimensional, 
then a unique branch of solutions exists in the vicinity of the
bifurcation. For example, in the problem of the buckling of a rod with
square cross-section, solutions exist in which the rod buckles
either parallel to two of the sides or diagonally. 
The equivariant branching lemma is known to generalize 
to the case where the dimension of the fixed-point subspace 
is odd, but it is in general untrue when the dimension is even. 

In this paper, the case of a two-dimensional fixed-point
subspace is considered, following on from earlier work.
Mcleod and Sattinger (1973) studied this problem, 
and showed that there are one or three solution branches for generic
values of the coefficients of the quadratic terms, corresponding to
five different phase portraits; however they regarded the double zero
eigenvalue as a degeneracy rather than a consequence of symmetry.
Lauterbach (1992), Leblanc {\it et al.} (1994) 
and Lari-Lavassani {\it et al.} (1994) 
formulated the problem in the context of equivariant bifurcation 
theory and obtained specific bifurcation theorems. 
In section~2 below, the basic theory of the subject 
is recalled, together with a summary of some of the known bifurcation
theorems. Section~3  gives specific conditions under which bifurcation
can be guaranteed, in terms of the coefficients of the quadratic
terms, and also gives the relationship between these coefficients and
the five possible phase portraits.
In section~4 these results are applied to the problem of bifurcation
from spherical symmetry, $O(3)$, in the fixed point subspace of the symmetry
group of the cube, which is two-dimensional for six of the irreducible
representations of $O(3)$.

\section{Branching in two dimensions}
\label{sec:2}

There is no simple generalization of the equivariant branching lemma
(EBL) to the case of a two-dimensional fixed point subspace.
In this section, some of the key definitions of the subject of
bifurcation with symmetry are recalled and some of the known results
are summarised. 

Near a multiple bifurcation point, 
the eigenfunctions form the basis of a vector space $V$. 
The symmetry group $G$ acts on  $V$  by multiplying 
points in $V$ by matrices that form a representation of $G$.
In most cases this representation is irreducible, 
meaning that there is no proper $G$-invariant subspace of $V$. 
The symmetry of a point $v$ in $V$ is described by its isotropy 
subgroup, defined by $\{g \in G: g v = v \}$. 

For any subgroup $H \subset G$, the fixed-point subspace of $H$ is
\begin{equation}
\mbox{Fix}(H) = \{ v \in V : h v = v \mbox{ for all } h \in H\}. 
\end{equation}
Fixed-point subspaces are important because they are invariant subspaces:
if an initial condition lies in $\mbox{Fix}(H)$ then the system
remains in $\mbox{Fix}(H)$ for all time. 

The normalizer $N(H)$ is defined by 
\begin{equation}
N(H) = \{g \in G : g^{-1} H g = H\}.
\end{equation}
This is the largest subgroup of $G$ in which $H$ is normal. 
Equivalently, this is the subgroup of $G$ which maps $\mbox{Fix}(H)$
to itself:
\begin{equation}
N(H) = \{g\in G : g y \in \mbox{Fix}(H) 
            \mbox{ for all } y \in \mbox{Fix}(H)\}.\label{neq}
\end{equation}
Clearly $H \subset N(H)$, and if  $N(H)=H$ then there is no symmetry
within $\mbox{Fix}(H)$, since all elements of $H$ act as the identity
in  $\mbox{Fix}(H)$. In general, the symmetry group acting on
$\mbox{Fix}(H)$ is the quotient group $N(H)/H$. 

Henceforth it will be assumed that the bifurcation problem with a
symmetry group $G$ satisfies the following conditions.

\begin{enumerate}
\item $G$ is a compact Lie group acting absolutely irreducibly on a
vector space $V$.
This implies that  $\mbox{Fix}(G) = 0$ or $\mbox{Fix}(G) = V$.
The latter case only holds for the trivial one-dimensional 
representation of $G$, so we may assume that  $\mbox{Fix}(G) = 0$.
\item The function $f: R \times V \rightarrow V$ is smooth and 
$G$-equivariant, i.e.\  $f(\lambda, g v) = g f(\lambda,v)$ for all $g
\in G$, $v\in V$. This implies that $f(\lambda, 0) = g f(\lambda,0)$,
and hence $f(\lambda, 0) \in \mbox{Fix}(G)$, so $f(\lambda,0) = 0$. 
Since the action of $G$ is absolutely irreducible, the linearization
of $f$ is a scalar multiple of the identity, $h(\lambda) I$.
\item $h(0)=0$, and $h'(0)\ne 0 $, i.e.\ there is a stationary
bifurcation at $\lambda=0$ and the eigenvalue passes transversely
through zero as $\lambda$ passes through zero.
\item $H$ is an isotropy subgroup of $G$ with $\mbox{Dim(Fix}(H)) = 2$. 
\item The quadratic terms in the Taylor expansion of $f$ do not all
vanish in $\mbox{Fix}(H)$.
\end{enumerate}

Note that the first three conditions are the standard ones for the
application of the EBL (Ihrig and Golubitsky 1984, Golubitsky,
Stewart and Schaeffer 1988, Chossat and Lauterbach 2000).

The equations $f(\lambda, v) = 0$ within $\mbox{Fix}(H)$ can then be
scaled, for small $\lambda$, by $\lambda = \epsilon \mu/h'(0)$, $v =
\epsilon x$, to give, as $\epsilon \rightarrow 0$, 
\begin{eqnarray}
0 & = & \mu x + a x^2 + b x y + c y^2, \label{e1}\\
0 & = & \mu y + p x^2 + q x y + r y^2,  \label{e2}
\end{eqnarray}
involving six constants $a$, $b$, $c$, $p$, $q$, $r$;
from the last of the above conditions,  not all of these are zero. 

These two equations were  discussed by Mcleod and
Sattinger (1973), who showed that for generic values 
of the six constants there exist either one or three 
solutions. However, a general bifurcation theorem cannot be deduced
from this, because in the context of equivariant bifurcation 
theory  the six constants are not independent.
There are two types of constraints on these coefficients
(see, for example, Stewart and  Dias 2000).
The first type of constraint arises from normaliser symmetries which
act on $\mbox{Fix}(H)$ as shown in (\ref{neq}); 
for example, there may be a normaliser symmetry acting as a reflection
$y\to -y$, in which case $p=b=r=0$. Secondly, there may be other
constraints arising from equivariance with respect to the original
symmetry group $G$; these are sometimes referred to as `hidden' symmetries.

An example to illustrate this point is the symmetry group of the
tetrahedron including reflections, which is isomorphic to the symmetric
group $S_4$ and also the group $O$ of rotations of the cube. 
In the natural three-dimensional irreducible representation, 
the group is generated by permutations of the three Cartesian coordinate
axes and a sign change of any two, and the bifurcation equations
consistent with these symmetries to quadratic order are 
\[
0 = \mu x + y z = \mu y + x z = \mu z + x y . 
\]
Note that there is only one equivariant quadratic term;
this is a feature common to many bifurcation problems,
including the cases of $S_N$ symmetry and $O(3)$ symmetry,
although there are also examples where this is not the case, such as
icosahedral symmetry (Hoyle 2003). 
The fixed point subspace of the symmetry group $Z_2$ generated by 
$y \leftrightarrow z$ 
is two-dimensional and the bifurcation equations reduce to 
\begin{equation}
0 = \mu x + y^2 = \mu y + x y. \label{tet2}
\end{equation}
The normaliser is $Z_2 \times Z_2$, including the element $(y,z)\to
(-y,-z)$, and so the normaliser symmetries force $p=b=r=0$ in
(\ref{tet2}), but we also have  
$a=0$ and $c=q$, arising from the original symmetry group.
Furthermore, there are only two non-zero solutions to (\ref{tet2}),
$x=-\mu$, $y=\pm \mu$, not one or three as expected from 
(\ref{e1}, \ref{e2}).
(In fact, in this example, a third solution in which $y=z=0$ does
exist when cubic terms are included, as required by the EBL since
this subspace is one-dimensional.)

Some rigorous results on bifurcation in two-dimensional subspaces are
as follows. A result of Lauterbach (1992) (Theorem 1.2),
obtained using degree theory, is that at least one 
non-zero solution exists if,
in addition to the five conditions above, 
\begin{equation}
Q^{-1}(0) = 0, \quad \mbox{i.e.} \quad Q = 0 \implies x=y=0, 
\label{laut}
\end{equation}
where $Q$ represents the quadratic terms in  (\ref{e1}, \ref{e2}).

Lauterbach's condition (\ref{laut}) is sufficient 
but not necessary. 
Leblanc {\it et al.} (1994) showed (Theorem 2.1) that in the case of
gradient dynamics, when $Q$ can be written as the gradient of a scalar 
quantity, a branch of solutions to the truncated problem
(\ref{e1}, \ref{e2}) always exists. Furthermore (Theorem 3.2),
they showed that a smooth branch of solutions of the non-truncated
problem exists if either $Q^{-1}(0) \ne 0$ or a fourth-order
polynomial in the six coefficients $a\ldots r$ is non-zero. 

Hence, by combining the results of Lauterbach (1992) 
and Leblanc {\it et al.} (1994), bifurcation to non-zero 
solutions always occurs occurs in the case of
gradient dynamics, since both the cases $Q^{-1}(0) = 0$ and 
$Q^{-1}(0) \ne 0$ are covered.

\section{Phase portraits}

In this section, solutions to the system (\ref{e1}, \ref{e2})
are considered in more detail.
First, the condition under which there are one or three roots is
derived.

Multiplying (\ref{e1}) by $y$, (\ref{e2}) by $x$ and subtracting 
yields the cubic equation 
\begin{equation}
0 = c z^3 + (b-r)z^2 + (a-q) z - p \label{cubic}
\end{equation}
for $z=y/x$, which generically has one or three
solutions. 
The condition for the existence of three roots depends on the quantity
\begin{equation}
P \equiv 27 c^2 p^2-18 c p (b-r)(q-a)-4 p (b-r)^3-4 c (q-a)^3-(q-a)^2
  (b-r)^2 . 
\label{3roots}
\end{equation}
There are three roots  if $P < 0$, one root if $P > 0$,
and a double root if $P = 0$. 
An alternative  derivation of the existence of one or three solutions,
employing the implicit function theorem, is given by Mcleod and
Sattinger (1973).  
Given a solution for $z = y/x$, substitution into 
(\ref{e1}) gives either the trivial solution $x=0$
or $x = -\mu/(a+bz+cz^2)$. This fails to give a finite solution if the
denominator is zero, and it is straightforward to show that this
occurs if the quantity 
\begin{equation}
R \equiv (a r-p c)^2 + (c q-b r) (a q-b p)  \label{res}
\end{equation}
is zero.  $R$ is the resultant of the two quadratic polynomials
appearing in (\ref{e1}, \ref{e2}). It is a well known result of
algebraic geometry that these two homogeneous polynomials have a
non-zero solution if and only if the resultant is zero (see, for
example, Cox, Little and O'Shea, 1998). 
Thus  the non-degeneracy condition $R\ne 0$ is exactly 
the condition (\ref{laut}), $Q^{-1}(0) = 0$. 
The degenerate case $R = 0$ corresponds to a solution at infinity. 

The non-degeneracy conditions giving one or three solutions 
are therefore $P\ne 0$ and $R\ne 0$. 
These conditions have not been given explicitly in earlier work,
although Leblanc {\it et al.} (1994) give $P$ in the special case of
gradient dynamics.

It is also of interest to consider the two-dimensional phase portraits
within $\mbox{Fix}(H)$, by considering the time dependent form of 
 (\ref{e1},  \ref{e2}):
\begin{eqnarray}
\dot x  & = & \mu x + a x^2 + b x y + c y^2, \label{et1}\\
\dot y  & = & \mu y + p x^2 + q x y + r y^2  \label{et2}. 
\end{eqnarray}
In this system, the straight line passing through the origin and any
fixed point is an invariant line, since if $z=y/x$ satisfies 
(\ref{cubic}) then $\dot z = (\dot y x - \dot x y)/x^2$
which is zero, since this is proportional to the same cubic appearing in
(\ref{cubic}). This observation forces all fixed points to be saddles
or nodes, and prohibits the existence of any periodic orbits. 
The invariance of these lines is a consequence of the truncation,
and is broken at higher order.

Stability of the fixed points is determined by the eigenvalues of the Jacobian
\begin{equation}
J =\pmatrix{{\mu + 2 a x + b y} & b x + 2 c y \cr
      2 p x + q y   & \mu + 2 r y + q x} . \label{jac}
\end{equation}
One eigenvalue is $-\mu$, corresponding to the eigenvector
lying along the invariant line.

In general, there are five possible different types of phase portrait
(McLeod and Sattinger 1973). 
Which of these occurs is determined by the values of 
$P$ and $R$ and a third quantity,
\begin{equation}
I = a q + b r - b p - c q .
\end{equation}
The following result gives the conditions on $P$, $R$ and $I$
that determine the phase portrait.

{\bf Theorem}\\
In the system (\ref{et1}, \ref{et2}),\\
(a) If $P > 0$ and $R > 0$ there is one fixed point, which is a
saddle point.\\
(b) If $P > 0$ and $R < 0$ there is one fixed point, which is a node.\\
(c) If $P < 0$ and $R > 0$ there are three fixed points, one node and
two saddles.\\
(d) If $P < 0$, $R < 0$ and $I > 0$ there are  three fixed points,
two  nodes and one saddle.\\
(e) If $P < 0$, $R < 0$ and $I < 0$ there are  three fixed points,
all of which are saddles.

{\bf Proof}\\
Note first that $P$, $R$ and $I$ are invariant under rotation of
coordinate axes. This can be verified by determining that under
an infinitesimal rotation through an angle $\theta$,
\[
\frac{d a}{d\theta} = b + p , \quad \frac{d b}{d\theta} = q+2c-2a, \quad 
\frac{d c}{d\theta} = r - b , 
\]
\[
\frac{d p}{d\theta} = q - a , \quad \frac{d q}{d\theta} = 2r-2p-b , \quad 
\frac{d r}{d\theta} = -c-q .
\]
To check that $I$ is invariant, we expand
\[
\frac{d I}{d\theta} = q (b+p)+(r-p)(q+2c-2a)-q (r - b)-b(q - a)
         +(a-c)(2r-2p-b ) - b(c+q)
\]
and verify that this quantity is zero, 
and similarly for $P$ and $R$.

Exploiting the rotation invariance and the fact that there is always 
at least one non-zero fixed point for $P\ne 0$, $R\ne 0$, 
we can always rotate coordinates so that a fixed point lies on the $x$
axis with $x < 0$. After this rotation, it follows from 
(\ref{et2}) that $p = 0$ and from (\ref{et1}) that the fixed point is
at $x=-\mu/a$. 
It can be assumed that $\mu > 0$, because the case $\mu < 0$
is equivalent under a sign change of $x$, $y$ and $t$,
and so $a > 0$. 
The eigenvalues of the fixed point at $x=-\mu/a$, $y=0$
are, from (\ref{jac}), $-\mu $ and $\mu (1-q/a)$,
so the fixed point is a saddle point if $a>q$ and a node if $a<q$.
Consider now the five cases above in turn.

(a) If $P > 0$ and $R > 0$,
then after rotating so that $p=0$ and $a>0$, $P>0$ gives
\begin{equation}
4c (a-q ) > (b-r)^2 \label{p01}
\end{equation}
so that either $c>0$ and $a>q$,  or $c<0$ and $a<q$.
The condition $R>0$ becomes
\begin{equation}
a r^2 + c q^2 - brq >0,
\end{equation}
which can be written as 
\[
(a-q)r^2 + c q^2 - rq(b-r) > 0. 
\]
Now if $c<0$ and $a<q$, this quadratic function of  $r$ must have 
real roots for the inequality to be satisfied, 
so $(b-r)^2 > 4 (a-q) c$, which contradicts (\ref{p01}).
Thus the only possibility is $c>0$ and $a>q$,
so the fixed point is a saddle point. 
The phase portrait in this case is shown in figure~1(a). 

(b) If $P > 0$ and $R < 0$ then (\ref{p01}) holds and
$R<0$ can be written as 
\[
(a-q)r^2 + c q^2 - rq(b-r) <  0. 
\]
Using the same argument as in the previous case,
if $a>q$ and $c > 0$
the inequality can only be satisfied 
if $(b-r)^2 > 4 (a-q) c$, contradicting (\ref{p01}).
Therefore $c<0$ and $a<q$, 
and the fixed point is a node, as shown in figure~1(b).

(c) If $P < 0$ and $R > 0$ there are three fixed points.
Using the result of Mcleod and Sattinger (1973) that at least 
one of these is a saddle point, a rotation of coordinates can be
carried out so that $p=0$ and a saddle point lies
at  $x=-\mu/a$, $y=0$; hence $a>q$. 
The quadratic equations for the coordinates $(x_1, y_1)$ 
and $(x_2, y_2)$ at the other fixed points, obtained from (\ref{e1})
and  (\ref{e2}), are
\begin{eqnarray}
(a r^2 + c q^2 - brq)x^2 + (r^2-br+2cq)\mu x + c\mu^2 & = & 0 \label{x2eq},\\
(a r^2 + c q^2 - brq)y^2 + (2ar-rq-bq)\mu y + (a-q)\mu^2  & = & 0 \label{y2eq}.
\end{eqnarray}
Since $R>0$, the leading term in each quadratic is positive,
and so from (\ref{y2eq}), $y_1 y_2 > 0$.
Thus all three fixed points lie in the same half plane, $y\geq 0$ or
$y \leq 0$. In the case that all lie in the upper half plane,
consider an anticlockwise rotation that brings the next fixed point
onto the negative axis. After this rotation, the two points not on the
$x$ axis have  $y_1 y_2 < 0$, so $a<q$ and so the point on the axis 
must be a node. After a further rotation so that the third fixed point
lies on the $x$ axis, both other fixed points have $y < 0$,
so $a>q$ and the point on the axis is a saddle. 
Hence for $P < 0$ and $R > 0$, there are three fixed point in the same
half plane, and the middle of the three is a node while 
the other two are saddle points (figure~1(c)).

(d) If $P < 0$, $R < 0$ and $I > 0$,
then $4c(a-q)<(b-r)^2$ and $a^2 r^2 + (cq - br) aq <0$,
so $cq-br$ and $aq$ must have opposite signs. 
Since $I>0$, $aq + br -cq > 0$ and so both $aq$ and $br -cq$ must be
positive. Since $P<0$, there are  three fixed points.
After a rotation of coordinates, one lies on the $x$ axis and the other
two satisfy (\ref{x2eq}, \ref{y2eq}). 
At the other two fixed points,  one eigenvalue is $-\mu$ and the other,
from (\ref{jac}), is $3\mu + (2a+q)x + (2r+b)y$. 
The sum of the `non-radial' eigenvalues at these fixed points
can be found using the sums of the roots obtained from (\ref{x2eq},
\ref{y2eq}).  The result is 
$q a \mu ( r^2 + b^2 -2br + 4 c q - 4a c ) /R$, 
which is negative. Thus, at least one of these fixed points 
must be a stable node. 
Given that one fixed point is a node, this point can be rotated onto
the $x$ axis and an argument exactly analogous to that of case 3 can
be followed, leading to the conclusions that the three fixed points lie
in the same half plane, the middle one is a saddle point and the other
two are nodes (figure~1(d)).

(e) If $P < 0$, $R < 0$ and $I < 0$ then, as in case (d), $cq-br$ and $aq$
must have opposite signs. In this case  $aq + br -cq < 0$,
so $aq< 0$ and $br - cq < 0$. Hence $a^2 - a q > 0$ and so 
the fixed point on the $x$ axis is a saddle point. 
Thus in this case all three fixed points are saddles,
as shown in figure~1(e).

{\bf Remarks}

\begin{enumerate}

\item 
The five possible phase portraits are illustrated in 
figure~1, for $\mu > 0$. 
The case $\mu < 0$ is obtained by reversing all the arrows 
and rotating through $\pi$. 

\item
Stable fixed points can only occur in cases (b), (c) and (d). 
Of course, these stability assignments only apply within 
Fix$(H)$, and there are directions transverse to Fix$(H)$ 
which may be stable or unstable.
In fact, it is generally the case that all solution branches 
obtained are unstable in the full system, since it can be shown that 
one eigenvalue is $-\mu$ but the sum of the eigenvalues is $n\mu$,
where $n$ is the dimension of the representation
(see, for example, Chossat {\it et al.} 1990).

\item
As $P$ passes through zero there is a saddle-node
bifurcation, leading to a transition from case (a) to case (c) or from 
case (b) to case (d). 
As $R$ passes through zero, a fixed point passes through infinity,
giving transitions between cases (a) and (b),
(c) and (d) or (c) and (e). There can be no direct transition between
cases (d) and (e) because $I=0$ is incompatible with $R<0$.

\item 
In the case of gradient-like dynamics, there are additional
constraints $b = 2p$ and $q = 2c$. 
Under these conditions it can easily be shown that $P>0 \implies R>0$
and $R < 0 \implies I < 0$, so that there are only three possible phase
portraits, (a), (c) and (e).

\end{enumerate}

\begin{figure}
\epsfxsize0.9\hsize\epsffile{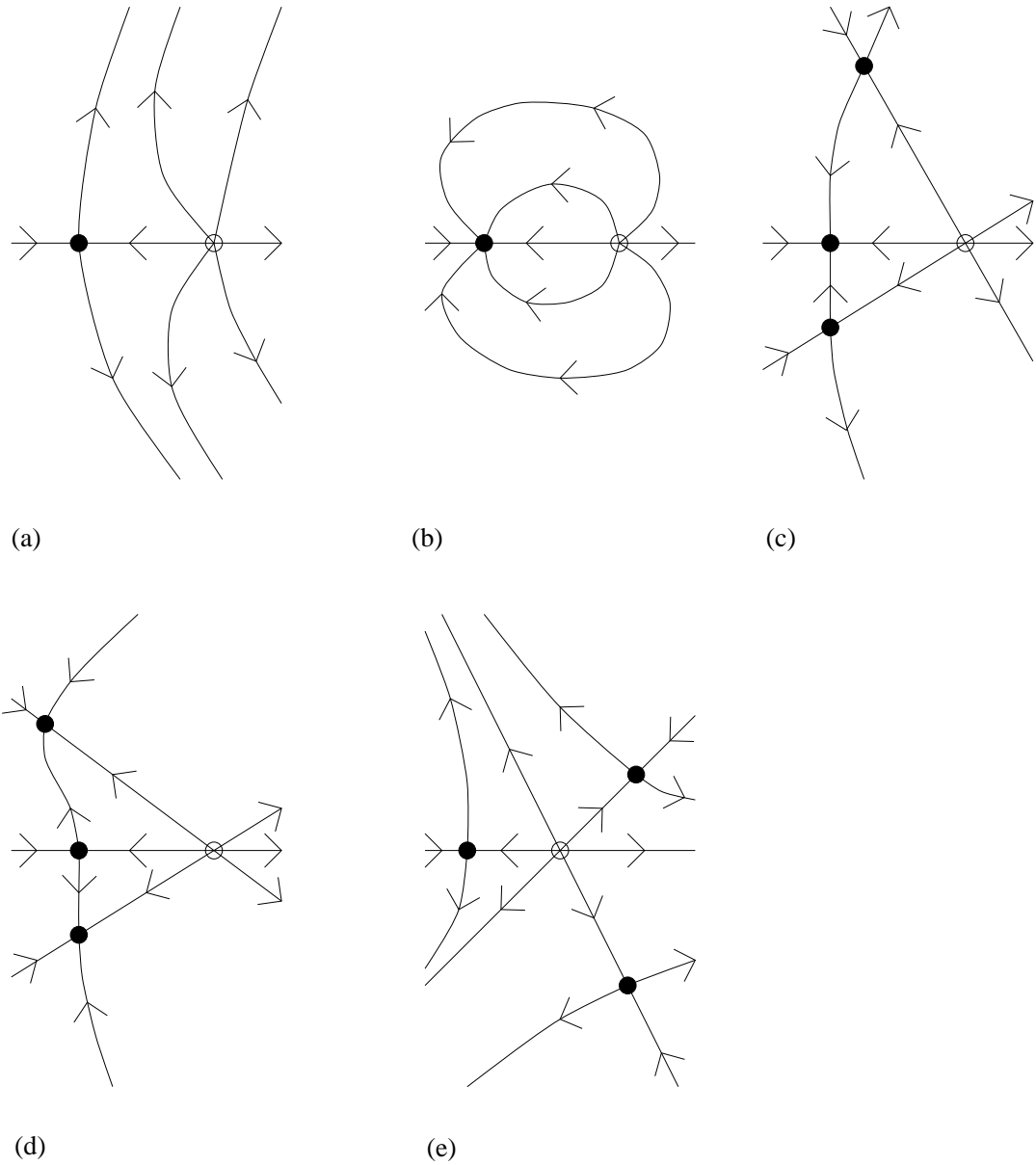}
\caption{The five possible types of phase portrait for 
(\ref{et1}, \ref{et2}), when $\mu>0$ (for $\mu < 0$, all arrows are reversed). 
The open circle  indicates the trivial solution and filled circles are 
non-zero solutions. The five case are (a) one solution, saddle,
(b) one solution, node, (c) three solutions, saddle, node, saddle,
(d) three solutions, node, saddle, node, (e) three solutions, all
saddles. }
\end{figure}

\section{Application to bifurcation with spherical symmetry}
\label{sec:app}

An interesting example of bifurcation with symmetry is the problem of
bifurcation from an initial state with spherical symmetry.
This problem has many physical applications and has been studied at
length by Sattinger (1978, 1979), Ihrig and Golubitsky (1984),
Chossat {\it et al.} (1990), and Matthews (2003). 
The irreducible representations have dimension $2l + 1$ and the
eigenfunctions are the spherical harmonics of degree $l$. 
Of particular interest are patterns with the symmetry groups of 
rotations of the Platonic solids: 
the tetrahedron, $T$, with 12 elements,
the cube / octahedron, $O$, with 24 elements, containing $T$,
and the icosahedron / dodecahedron, $I$, with 60 elements, also
containing $T$.  
When $l$ is odd, the quotient group $N(H)/H$ always includes 
the group $Z_2^C$, generated by point inversion through the origin, 
corresponding to a rotation through $\pi$ if  $\mbox{Dim(Fix}(H)) =2$,
so the quadratic terms are all forced to be zero by the 
normaliser symmetry.
For even $l$, the isotropy subgroups are all of the form
$H \oplus Z_2^C$, where $H$ is a subgroup of SO(3) and $\oplus$
indicates a direct product. 
The groups $O\oplus Z_2^C$ and  $I\oplus Z_2^C$ are maximal isotropy
subgroups, and satisfy the condition $N(H)=H$.
The dimensions of their fixed-point subspaces
are $[l/4] + [l/3] - l/2 + 1$ for $O\oplus Z_2^C$
and $[l/5] + [l/3] - l/2 + 1$ for $I\oplus Z_2^C$,
where $[x]$ is the largest integer less than or equal to $x$. 

The dynamics is gradient-like, so the following two existence results
follow directly from the observation at the end of 
section~2, by determining the values of $l$ for which
 $\mbox{Dim(Fix}(H)) = 2$. 
\begin{enumerate}
\item There exists at least one, but not more than  three branches of solutions with isotropy $O\oplus Z_2^C$ for $l=12$, 16, 18, 20, 22, 26. 
\item There exists at least one, but not more than  three  branches of
solutions with isotropy  $I\oplus Z_2^C$ for $l=30$, 36, 40, 42, 46,
48, 50, 52, 54, 56, 58, 62, 64, 68, 74.
\end{enumerate}

Since the dynamics is gradient-like, the phase portraits
within the two-dimensional subspace  $\mbox{Fix}(H)$
may have any of the three forms shown in figure~1(a), (c) or (e) 
provided that the non-degeneracy conditions $P\ne 0$, $R\ne 0$ are
satisfied. 
For the case $l=12$, $H = O\oplus Z_2^C$,
the computation of the coefficients was carried out by
Leblanc {\it et al.} (1994); the result is, after 
removing some common factors,
\[
a = - \frac{9913 \sqrt{8398}} {390}, \quad
c = \frac{73501\sqrt{8398}} {4940}, \quad
p = \frac{104\sqrt{5313}} {3}, \quad
r = -\frac{339119\sqrt{5313}} {4370},
\]
with $b=2p$ and $q=2c$. 
Calculating the values of $P$, $R$ and $I$ shows that all three are
negative. Hence there are three distinct branches of solutions with
isotropy $O\oplus Z_2^C$ for $l=12$, and within this subspace, all three are
saddle points. 
 
A similar calculation was carried out for the other values of $l$ for
which $\mbox{Dim(Fix}(O\oplus Z_2^C)) = 2$, using the
computer algebra package Maple to do the cumbersome manipulations. 
First, the $2l+1$ equations for the amplitudes of the spherical
harmonics $Y_l^m(\theta,\phi)$, $m=-l\ldots l$ are constructed using
the fact that the quadratic coefficients are the Clebsch--Gordan
coefficients (Sattinger, 1978). The subspace $\mbox{Fix}(D_4\oplus Z_2^C)$ 
of dimension $1 + [l/4]$ is found simply by restricting to those modes
in which $m$ is a multiple of four. The subspace $\mbox{Fix}(O\oplus
Z_2^C)$  is then obtained from this by imposing invariance under a 
rotation of $\pi/2$. 
The results of these calculations are as follows.
For $l=12$, $16$, $22$ and $26$,  $P$, $R$ and $I$
are all negative, so there are three stationary solutions
and each solution is a saddle point in $\mbox{Fix}(O\oplus Z_2^C)$,
as in figure~1(e).
For $l=18$, $P<0$ and $R>0$, so there are three solutions, one
of which is a node and the other two are saddles (figure~1(c)). 
For $l=20$, $P>0$ and $R>0$,  so there is only one solution, 
which is a saddle point.

\section*{References}

Chossat, P. and Lauterbach, R., 2000, {\it Methods in Equivariant
Bifurcations and Dynamical Systems} (World Scientific).

Chossat, P.,  Lauterbach, R.  and  Melbourne, I., 1990,
Steady-state bifurcation with O(3) symmetry.
{\em Arch. Rat. Mech. Anal.} {\bf 113}, 313--376.

Cicogna, G. 1981, Symmetry breakdown from bifurcation.
{\it Lettere al Nuovo Cimento} {\bf 31}, 600--602. 

Cox, D., Little, J. and O'Shea, D., 1998,
{\it Using Algebraic Geometry} (Springer).

Golubitsky, M., Stewart, I. and Schaeffer, D.G., 1988, {\it
Singularities and Groups in Bifurcation Theory} (Springer). 

Hoyle, R.B., 2003, Shapes and cycles arising at the steady bifurcation
with icosahedral symmetry.  {\it Physica} D, submitted.

Ihrig, E. and Golubitsky, M., 1984,
Pattern selection with O(3) symmetry.
{\it Physica} {\bf 13D}, 1--33. 

Lari-Lavassani, A., Langford, W.F. and Huseyin, K., 1994,
Symmetry-breaking bifurcations on multidimensional fixed point subspaces.
{\it Dynamics and Stability of Systems} {\bf 9}, 345--373.

Lauterbach, L., 1992,
Spontaneous symmetry breaking in higher dimensional fixed point spaces
{\it ZAMP} {\bf 43}, 430--448.

LeBlanc, V.G., Lari-Lavassani, A. and Langford, W.F., 1994,
Symmetry-breaking for leading order gradient maps in 
$R^2$ with applications to $O(3)$.
{\it Nonlinearity} {\bf 7}, 577--594.

Matthews, P.C., 2003,
Transcritical bifurcation with $O(3)$ symmetry.
{\it Nonlinearity}, in press.

McLeod, J.B. and Sattinger, D.H., 1973,
Loss of stability and bifurcation at a double eigenvalue.
{\it J. Functional Analysis}, {\bf 14}, 62--84.

Sattinger, D., 1978,
Bifurcation from rotationally invariant states.
{\it J. Math. Phys.} {\bf 19}, 1720--1732.

Sattinger, D., 1979, {\it Group Theoretic Methods in Bifurcation
Theory} (Springer).

Stewart, I. and  Dias, A.P.S., 2000,
Hilbert series for equivariant mappings restricted to invariant hyperplanes
{\it J. Pure and Applied Algebra} {\bf 151}, 89--106.

Vanderbauwhede, A., 1980, Symmetry and bifurcation near families of
solution. {\it J. Differential Equations} {\bf 36}, 173--178. 

\end{document}